\documentclass[11pt,a4paper]{article}
\usepackage{amsmath,amsthm,amssymb}
\usepackage[top=2.5cm,right=3.5cm,left=3.5cm,bottom=3.5cm]{geometry}

\title{A refinement of Ramanujan's factorial approximation}

\author{Michael D. Hirschhorn$^1$ and Mark B. Villarino$^2$\\[6pt]
$^1$School of Mathematics and Statistics, University of New South Wales,\\
Sydney NSW 2052, Australia\\[12pt]
$^2$Escuela de Matem\'atica, Universidad de Costa Rica\\
San Jos\'e 11501, Costa Rica}

\date{}

\newtheorem{thm}{Theorem}

\newtheorem{propn}[thm]{Proposition}

\newtheorem{defn}[thm]{Definition}

\newcommand{\bZ}{\mathbb{Z}}       
\newcommand{\dsp}{\displaystyle}   
\newcommand{\hideqed}{\renewcommand{\qed}{}} 
\newcommand{\word}[1]{\quad\text{#1}\quad} 

\newcommand{\half}{{\mathchoice{\thalf}{\thalf}{\shalf}{\shalf}}}
\newcommand{\shalf}{{\scriptstyle\frac{1}{2}}} 
\newcommand{\thalf}{{\tfrac{1}{2}}} 

\begin{document}

\maketitle

\section{Introduction}

All we have of Ramanujan's work in the last year of his life is about
100 pages (probably a small fraction of his final year's output), held
by Trinity College, Cambridge, and named by George E. Andrews
``Ramanujan's Lost Notebook''. It was published in photocopied
form~\cite{RamanLost}. In it, Ramanujan \cite[p.~339]{RamanLost} makes
the claim that
$$
\Gamma(x + 1) = \sqrt{\pi} \biggl( \frac{x}{e} \biggr)^x
\biggl( 8x^3 + 4x^2 + x + \frac{\theta_x}{30} \biggr)^\frac{1}{6}
$$
where $\theta_x \to 1$ as $x \to \infty$ and
$\dfrac{3}{10} < \theta_x < 1$, and gives some numerical evidence
for this last statement.

Inspired by this, we confine ourselves to the positive integers, and
prove the following stronger result.

\begin{thm} 
\label{th:main}
Let the function, $\theta(n)\equiv \theta_n$, be defined for $n = 1,2,\dots$ by the
equation:
$$
n! := \sqrt{\pi} \biggl( \frac{n}{e} \biggr)^n
\biggl(8n^3 + 4n^2 + n + \frac{\theta(n)}{30} \biggr)^\frac{1}{6}.
$$
Then, the correction term $\theta(n)$
\begin{enumerate}
\item 
satisfies the inequalities:
\begin{equation}
\label{eq:ineq} 
1 - \frac{11}{8n} + \frac{79}{112n^2} < \theta(n)
< 1 - \frac{11}{8n} + \frac{79}{112n^2} + \frac{20}{33n^3}\,;
\end{equation}
\item 
is an \textbf{increasing} function of $n$; and
\item 
is \textbf{concave}, that is,
$$
\theta_{n+1} - \theta_n < \theta_n - \theta_{n-1}.
$$
\end{enumerate}
\end{thm}

The inequalities \eqref{eq:ineq} are new. In 2006, Hirschhorn
\cite{HirschhornNewStir} proved a less exact version of the
inequalities \eqref{eq:ineq}. In 2001, Karatsuba~\cite{Karatsuba} proved
Ramanujan's approximation and gave a  proof, quite different from ours, of the
monotonicity of the correction term~$\theta_x$, for all real
$x\geq 1$, a result which is stronger than ours. Moreover, although Karatsuba derived an asymptotic
expansion for $\theta_x$, including a uniform error term, she did not
derive any explicit numerical inequalities, as we do. The monotonicity
of $\theta(n)$ was proved by Villarino, Campos-Salas, and
Carvajal-Rojas in~\cite{VillarinoCC} as a simple consequence of the
inequality in~\cite{HirschhornNewStir}; in that paper, the concavity
of $\theta(n)$ was also noted, without proof.

Our proofs use nothing more than the series for $\log{(1 + x)}$ and
$\exp\{x\}$.

We will find it convenient to use the following notation.
\begin{defn} The notation $$P_k(n)$$ means a  polynomial of degree $k$ in $n$ with all of its non-zero coefficients \textbf{\emph{positive}}.

\end{defn}

\section{The proofs}


\begin{propn} 
\label{pr:two}
The following inequality is valid for $n = 1,2,\dots$:
\begin{equation}
\label{eq:easy} 
0 < \biggl( n + \frac{1}{2} \biggr) 
\log \biggl( 1 + \frac{1}{n} \biggr) - 1
< \frac{1}{12n} - \frac{1}{12(n + 1)} \,.
\end{equation}
\end{propn}

\begin{proof}
We have, for $|x| < 1$,
$$
\log{(1 + x)} = x - \frac{x^2}{2} + \frac{x^3}{3} - \cdots \,.
$$
It follows that for $|x| < 1$,
$$
\log \biggl( \frac{1 + x}{1 - x} \biggr)
= 2 \biggl( x + \frac{x^3}{3} + \frac{x^5}{5} + \cdots \biggr).
$$
If we set $x = \dfrac{1}{2n + 1}$ where $n \in \bZ^+$, we obtain 
\begin{equation}
\label{eq:log-roll} 
\log\biggl( 1 + \frac{1}{n} \biggr)
= 2 \biggl( \frac{1}{2n + 1} + \frac{1}{3(2n + 1)^3}
+ \frac{1}{5(2n + 1)^5} + \cdots \biggr).
\end{equation}
It follows that
$$
\biggl( n + \frac{1}{2} \biggr) \log \biggl( 1 + \frac{1}{n} \biggr)
= 1 + \frac{1}{3(2n + 1)^2} + \frac{1}{5(2n + 1)^4} + \cdots
$$
Therefore
$$
0 < \biggl( n + \frac{1}{2} \biggr)
\log \biggl( 1 + \frac{1}{n} \biggr) - 1
< \frac{1}{3(2n + 1)^2} \cdot \frac{1}{1 - \dfrac{1}{(2n + 1)^2}}
= \frac{1}{12n} - \frac{1}{12(n + 1)} \,.
\eqno \qed
$$
\hideqed
\end{proof}

The inequality \eqref{eq:easy} leads to the following well-known version
of Stirling's formula.


\begin{propn} 
\label{pr:three}
The following inequality is valid for $n = 1,2,\dots$:
\begin{equation}
\label{eq:Stir-one} 
\sqrt{2\pi n} \biggl( \frac{n}{e} \biggr)^n < n!
\leq \sqrt{2\pi n} \biggl( \frac{n}{e} \biggr)^n
\exp \biggl\{ \frac{1}{12n} \biggr\} \,.
\end{equation}
\end{propn}

\begin{proof}
Let 
$$
a_n = n! \Big/ \sqrt{n} \biggl( \frac{n}{e} \biggr)^n.
$$
Then
$$
\frac{a_n}{a_{n + 1}}
= \biggl( 1 + \frac{1}{n} \biggr)^{n+\half} \! \bigg/ e
= \exp \biggl\{ \biggl( n + \frac{1}{2} \biggr) 
\log \biggl( 1 + \frac{1}{n} \biggr) - 1 \biggr\}.
$$
{}From \eqref{eq:easy} we have
\begin{equation}
\label{eq:ratio} 
1 < \frac{a_n}{a_{n+1}}
< \exp \biggl\{ \frac{1}{12n} - \frac{1}{12(n + 1)} \biggr\}.
\end{equation}
So $a_n$ is decreasing and, if we write $1$, $2$,\dots, $n - 1$
for~$n$ and multiply the results, we find
$$
\frac{a_1}{a_n} < \exp \biggl\{ \frac{1}{12}  -\frac{1}{12n} \biggr\}
< \exp \biggl\{ \frac{1}{12} \biggr\},
$$
or,
$$
a_n > a_1 \exp \biggl\{ -\frac{1}{12} \biggr\}
= \exp \biggl\{ \frac{11}{12} \biggr\}.
$$
It follows that $a_\infty = \lim_{n\to\infty} a_n$ exists, and
$$
a_\infty \geq \exp \biggl\{ \frac{11}{12} \biggr\}.
$$
In fact, from Wallis's product, $a_\infty = \sqrt{2\pi}$.

If in \eqref{eq:ratio} we write $n$, $n + 1$,\dots, $N - 1$ for~$n$,
multiply the results, let $N\to\infty$ and use Wallis's product, we
obtain, successively,
\begin{align*}
1 &< \frac{a_n}{a_N} 
< \exp \biggl\{ \frac{1}{12n} - \frac{1}{12N} \biggr\},
\\
a_N &< a_n < a_N \exp \biggl\{ \frac{1}{12n} - \frac{1}{12N} \biggr\},
\\
\sqrt{2\pi} &< a_n 
\leq \sqrt{2\pi} \exp \biggl\{ \frac{1}{12n} \biggr\}
\\
\intertext{and}
\sqrt{2\pi n} \biggl( \frac{n}{e} \biggr)^n 
&< n! \leq \sqrt{2\pi n} \biggl( \frac{n}{e} \biggr)^n
\exp \biggl\{ \frac{1}{12n} \biggr\}.
\tag*{\qed}
\end{align*}
\hideqed
\end{proof}


We shall improve on \eqref{eq:easy}, and so also
on~\eqref{eq:Stir-one}, by proving the next inequality.

\begin{propn} 
\label{pr:four}
The following inequality is valid for $n = 1,2,\dots$:
\begin{align}
& \frac{1}{12} \biggl( \frac{1}{n} -\frac{1}{n + 1} \biggr)
- \frac{1}{360} \biggl( \frac{1}{n^3} - \frac{1}{(n + 1)^3} \biggr)
+ \frac{1}{1260} \biggl( \frac{1}{n^5} - \frac{1}{(n + 1)^5} \biggr)
\nonumber \\
&\hspace*{20em}
- \frac{1}{1680} \biggl( \frac{1}{n^7} - \frac{1}{(n + 1)^7} \biggr)
\nonumber \\
&\qquad < \biggl( n + \frac{1}{2} \biggr)
\log \biggl( 1 + \frac{1}{n} \biggr) - 1
\nonumber \\[\jot]
&\qquad < \frac{1}{12} \biggl( \frac{1}{n} - \frac{1}{n + 1} \biggr)
- \frac{1}{360} \biggl( \frac{1}{n^3} - \frac{1}{(n + 1)^3} \biggr)
+ \frac{1}{1260} \biggl( \frac{1}{n^5} - \frac{1}{(n + 1)^5} \biggr).
\label{eq:log-bounds} 
\end{align}
\end{propn}

\begin{proof}
We have, from \eqref{eq:log-roll},
\begin{align*}
& \biggl(n + \frac{1}{2}\biggr) \log \biggl(1 + \frac{1}{n}\biggr) - 1
< \frac{1}{3(2n + 1)^2} + \frac{1}{5(2n + 1)^4}
+ \frac{1}{7(2n + 1)^6} \cdot \frac{1}{1 - \dfrac{1}{(2n + 1)^2}}
\\
&\qquad = \frac{1}{3(2n + 1)^2} + \frac{1}{5(2n + 1)^4}
+ \frac{1}{28n(n + 1)(2n + 1)^4}
\\[\jot]
&\qquad = \frac{1}{12} \biggl( \frac{1}{n} - \frac{1}{n + 1} \biggr)
- \frac{1}{360} \biggl( \frac{1}{n^3} - \frac{1}{(n + 1)^3} \biggr)
+ \frac{1}{1260} \biggl( \frac{1}{n^5} - \frac{1}{(n + 1)^5} \biggr)
\\[\jot]
&\qquad-\frac{163n^6 + 489n^5 + 604n^4 + 393n^3 + 141n^2 + 26n + 2}
{2520n^5 (n + 1)^5 (2n + 1)^4}\\
&\qquad < \frac{1}{12} \biggl( \frac{1}{n} - \frac{1}{n + 1} \biggr)
- \frac{1}{360} \biggl( \frac{1}{n^3} - \frac{1}{(n + 1)^3} \biggr)
+ \frac{1}{1260} \biggl( \frac{1}{n^5} - \frac{1}{(n + 1)^5} \biggr)
\end{align*}
and
\begin{align*}
& \biggl(n + \frac{1}{2}\biggr) \log\biggl(1 + \frac{1}{n}\biggr) - 1
> \frac{1}{3(2n + 1)^2} + \frac{1}{5(2n + 1)^4} + \frac{1}{7(2n + 1)^6}
+ \frac{1}{9(2n + 1)^8}
\\
&\qquad =\frac{1}{12} \biggl( \frac{1}{n} -\frac{1}{n + 1} \biggr)
- \frac{1}{360} \biggl( \frac{1}{n^3} - \frac{1}{(n + 1)^3} \biggr)
+ \frac{1}{1260} \biggl( \frac{1}{n^5} - \frac{1}{(n + 1)^5} \biggr)
\\
&\hspace*{4em}
- \frac{1}{1680} \biggl( \frac{1}{n^7} - \frac{1}{(n + 1)^7} \biggr)\\
&\hspace*{4em}
+\frac{P_{12}(n)}{5040n^7(n+1)^7(2n+1)^8}
\\
&\qquad > \frac{1}{12} \biggl( \frac{1}{n} -\frac{1}{n + 1} \biggr)
- \frac{1}{360} \biggl( \frac{1}{n^3} - \frac{1}{(n + 1)^3} \biggr)
+ \frac{1}{1260} \biggl( \frac{1}{n^5} - \frac{1}{(n + 1)^5} \biggr)
\\
&\hspace*{4em}
- \frac{1}{1680} \biggl( \frac{1}{n^7} - \frac{1}{(n + 1)^7} \biggr).
\end{align*}This completes the proof.
\end{proof}


We now demonstrate the greatly improved version of Stirling's formula.

\begin{propn} 
\label{pr:five}
For $n = 1,2,\dots$, the following inequality is valid:
\begin{align}
& \sqrt{2\pi n} \biggl( \frac{n}{e} \biggr)^n
\exp \biggl\{ \frac{1}{12n} - \frac{1}{360n^3} + \frac{1}{1260n^5}
- \frac{1}{1680n^7} \biggr\}
\nonumber \\[\jot]
&\qquad \leq n!
\leq \sqrt{2\pi n} \biggl( \frac{n}{e} \biggr)^n
\exp \biggl\{ \frac{1}{12n} - \frac{1}{360n^3}
+ \frac{1}{1260n^5} \biggr\}.
\label{eq:Stir-two} 
\end{align}
\end{propn}

\begin{proof}
It follows from \eqref{eq:log-bounds} that
\begin{align*}
& \exp\biggl\{ \frac{1}{12} \biggl(\frac{1}{n} -\frac{1}{n + 1}\biggr)
- \frac{1}{360} \biggl( \frac{1}{n^3} - \frac{1}{(n + 1)^3} \biggr)
+ \frac{1}{1260} \biggl( \frac{1}{n^5} - \frac{1}{(n + 1)^5} \biggr)
\\
&\hspace*{22em} - \frac{1}{1680}
\biggl( \frac{1}{n^7} - \frac{1}{(n + 1)^7} \biggr) \biggr\}
\\
&\qquad < \frac{a_n}{a_{n+1}}
\\
&\qquad < \exp \biggl\{ \frac{1}{12} 
\biggl( \frac{1}{n} -\frac{1}{n + 1} \biggr)
- \frac{1}{360} \biggl( \frac{1}{n^3} - \frac{1}{(n + 1)^3} \biggr)
+ \frac{1}{1260} \biggl( \frac{1}{n^5} - \frac{1}{(n + 1)^5} \biggr)
\biggr\}.
\end{align*}
Thus, for $N > n$,
\begin{align*}
& \exp \biggl\{ \frac{1}{12} \biggl( \frac{1}{n} -\frac{1}{N} \biggr)
- \frac{1}{360} \biggl( \frac{1}{n^3} - \frac{1}{N^3} \biggr)
+ \frac{1}{1260} \biggl( \frac{1}{n^5} - \frac{1}{N^5} \biggr)
- \frac{1}{1680} \biggl( \frac{1}{n^7} - \frac{1}{N^7} \biggr) \biggr\}
\\[\jot]
&\qquad < \frac{a_n}{a_N}
< \exp \biggl\{
\frac{1}{12} \biggl( \frac{1}{n} -\frac{1}{N} \biggr)
- \frac{1}{360}\biggl( \frac{1}{n^3} - \frac{1}{N^3} \biggr)
+ \frac{1}{1260} \biggl( \frac{1}{n^5} - \frac{1}{N^5} \biggr) \biggr\},
\end{align*}
and
\begin{align*}
& \sqrt{2\pi n} \biggl( \frac{n}{e} \biggr)^n 
\exp \biggl\{ \frac{1}{12n} - \frac{1}{360n^3} + \frac{1}{1260n^5}
- \frac{1}{1680n^7} \biggr\}
\\[\jot]
&\qquad \leq n!
\leq \sqrt{2\pi n} \biggl( \frac{n}{e} \biggr)^n \exp \biggl\{
\frac{1}{12n} - \frac{1}{360n^3} + \frac{1}{1260n^5} \biggr\}.
\tag*{\qed}
\end{align*}
\hideqed
\end{proof}


Extracting the fraction $\smash{\dfrac{1}{6}}$ from the exponents, we
see that we can write this last inequality in the form
\begin{align*}
& \sqrt{2\pi n} \biggl( \frac{n}{e} \biggr)^n
\biggl( \exp\biggl\{ \frac{1}{2n} - \frac{1}{60n^3} + \frac{1}{210n^5}
- \frac{1}{280n^7} \biggr\} \biggr)^\frac{1}{6}
\\[\jot]
&\qquad \leq n!
\leq \sqrt{2\pi n} \biggl( \frac{n}{e} \biggr)^n
\biggl( \exp\biggl\{ \frac{1}{2n} - \frac{1}{60n^3}
+ \frac{1}{210n^5} \biggr\} \biggr)^\frac{1}{6}.
\end{align*}


We now obtain upper and lower bounds for these new exponents.

\begin{propn} 
\label{pr:six}
The following inequalities are valid for  $n \geq 2$:
\begin{align}
& \exp\biggl\{ \frac{1}{2n} - \frac{1}{60n^3} + \frac{1}{210n^5}
- \frac{1}{280n^7} \biggr\}
\nonumber \\
&\hspace*{4em}
> 1 + \frac{1}{2n} + \frac{1}{8n^2} + \frac{1}{240n^3}
- \frac{11}{1920n^4} + \frac{79}{26880n^5}
\label{eq:lower-bound} 
\\
\intertext{and}
& \exp\biggl\{ \frac{1}{2n} - \frac{1}{60n^3} + \frac{1}{210n^5}
\biggr\}
\nonumber \\
&\hspace*{4em}
< 1 + \frac{1}{2n} + \frac{1}{8n^2} + \frac{1}{240n^3}
- \frac{11}{1920n^4} + \frac{79}{26880n^5} + \frac{1}{396n^6} \,.
\label{eq:upper-bound} 
\end{align}
\end{propn}


Assuming for the moment that these bounds are valid, we can now prove
the main result of this paper.

\begin{proof}[Proof of Theorem \ref{th:main}]
It follows from Proposition~\ref{pr:six} that for $n \geq 2$,
\begin{align*}
& \sqrt{2\pi n} \biggl( \frac{n}{e} \biggr)^n
\biggl( 1 + \frac{1}{2n} + \frac{1}{8n^2} + \frac{1}{240n^3}
- \frac{11}{1920n^4} + \frac{79}{26880n^5} \biggr)^\frac{1}{6}
\\
&\quad < n!
< \sqrt{2\pi n} \biggl( \frac{n}{e} \biggr)^n
\biggl( 1 + \frac{1}{2n} + \frac{1}{8n^2} + \frac{1}{240n^3}
- \frac{11}{1920n^4} + \frac{79}{26880n^5} + \frac{1}{396n^6}
\biggr)^\frac{1}{6},
\end{align*}
or:
\begin{align*}
& \sqrt{\pi} \biggl( \frac{n}{e} \biggr)^n
\biggl( 8n^3 + 4n^2 + n + \frac{1}{30} \biggl(
1 - \frac{11}{8n} + \frac{79}{112n^2} \biggr) \biggr)^\frac{1}{6}
\\
&\qquad < n!
< \sqrt{\pi} \biggl( \frac{n}{e} \biggr)^n
\biggl( 8n^3 + 4n^2 + n + \frac{1}{30} \biggl(
1 - \frac{11}{8n} + \frac{79}{112n^2} + \frac{20}{33n^3} \biggr)
\biggr)^\frac{1}{6}.
\end{align*}
This beautiful formula is the refined estimate~\eqref{eq:ineq}. It is
easy to check this for $n = 1$ also, so we have the desired result.

To show that $\theta(n)$ is increasing, from~\eqref{eq:ineq} it follows
that
\begin{align*}
\theta_{n+1} - \theta_n 
& > 1 - \frac{11}{8(n + 1)} + \frac{79}{112(n + 1)^2} - \biggl(
1 - \frac{11}{8n} + \frac{79}{112n^2} + \frac{20}{33n^3} \biggr)
\\[\jot]
&= \frac{(5082n^2 + 7792n + 8497)(n - 2) + 14754}{3696n^3(n + 1)^2}
\\[\jot]
&> 0 \word{for} n \geq 2,
\end{align*}
and it is easily checked for $n = 1$ also, so $\theta_n$ is increasing.

\vspace{6pt}

Finally, to prove the concavity of $\theta(n)$, we note that:
\begin{align*}
& \theta_{n+1} - 2\theta_n + \theta_{n-1}
\\
&\qquad < 1 - \frac{11}{8(n + 1)} + \frac{79}{112(n + 1)^2}
+ \frac{20}{33(n + 1)^3}
\\
&\hspace*{4em}
+ 1 - \frac{11}{8(n-1)} + \frac{79}{112(n-1)^2} + \frac{20}{33(n-1)^3}
- 2 \biggl( 1 - \frac{11}{8n} + \frac{79}{112n^2} \biggr)
\\[\jot]
&\qquad = - \frac{(2842n^4 + 6389n^3 + 15061n^2 + 85733n + 433747)
(n - 5) + 2166128}{1848n^2 (n - 1)^3 (n + 1)^3}
\\
&\qquad < 0 \word{for}  n \geq 5,
\end{align*}
and is easily checked for $n = 2$, $3$ and~$4$ also. 
\end{proof}


We complete the proof of the exponential inequalities as follows.

\begin{proof}[Proof of Proposition \ref{pr:six}]
Let $\dsp q := \frac{1}{2n} - \frac{1}{60n^3} + \frac{1}{210n^5}
- \frac{1}{280n^7}\,$. Then $q > 0$, and
\begin{align*}
\exp\{q\}
&> 1 + \frac{q}{1!} + \frac{q^2}{2!} + \frac{q^3}{3!} + \frac{q^4}{4!}
+ \frac{q^5}{5!}
\\
&=1 + \frac{1}{2n} + \frac{1}{8n^2} + \frac{1}{240n^3}
- \frac{11}{1920n^4} + \frac{79}{26880n^5} \\
&+\frac{P_{28}(n)(n-2)+ 5421638789368547485949} {50185433088000000 n^{35}}\\
&>1 + \frac{1}{2n} + \frac{1}{8n^2} + \frac{1}{240n^3}
- \frac{11}{1920n^4} + \frac{79}{26880n^5}
\end{align*}
which proves \eqref{eq:lower-bound} for $n\geq 2$.
Now let
$\dsp r := \frac{1}{2n} - \frac{1}{60n^3} + \frac{1}{210n^5}\,$. Then
$r > 0$, and
\begin{align*}
\exp\{r\}
&< 1 + \frac{r}{1!} + \frac{r^2}{2!} + \frac{r^3}{3!} + \frac{r^4}{4!}
+ \frac{r^5}{5!} + \frac{r^6}{6!} + \frac{r^7}{6!} +\cdots
\\[\jot]
&= 1 + \frac{r}{1!} + \frac{r^2}{2!} + \frac{r^3}{3!} + \frac{r^4}{4!}
+ \frac{r^5}{5!} + \frac{r^6}{6!} \bigg/ (1 - r)
\\
&=1 + \frac{1}{2n} + \frac{1}{8n^2} + \frac{1}{240n^3}
- \frac{11}{1920n^4} + \frac{79}{26880n^5} + \frac{1}{396n^6}\\
&-\frac{P_{23}(n)(n-3)+ 239259521624400145687307843} {20701491148800000 n^{25}(420n^5-210n^4+7n^2-2)}\\
&< 1 + \frac{1}{2n} + \frac{1}{8n^2} + \frac{1}{240n^3}
- \frac{11}{1920n^4} + \frac{79}{26880n^5} + \frac{1}{396n^6}
\word{for}  n \geq 3,
\end{align*}
which proves \eqref{eq:upper-bound} for $n \geq 3$. The case $n = 2$ is easily checked.
\end{proof}

\section{Final Remarks}

 \emph{Proposition 4} and \emph{Proposition 6} are special cases of the general expansion of $n!$, with an error term, which can be proved by using the Euler--Maclaurin sum formula. However, our proofs
are much more elementary, and can be extended to any degree of
accuracy desired. Still, our proofs do not supply the general formula
for the coefficients in the exponential version, although perhaps they
can be properly modified to do so.

Also, our technique for proving the positivity of certain large degree
polynomials seems to argue for a general property of polynomials
$P(x)$ with real coefficients that are positive for $x \geq a$. The
property in question is that there exists a $b \geq a$ such that the
quotient polynomial $Q(x)$ in the division algorithm
$P(x) \equiv Q(x)(x - b) + R$ has all its coefficients positive.

Finally, the inequalities \eqref{eq:ineq} can, with more work, be
extended to degrees three, four, etc., where the main coefficients are
given by the formula of Karatsuba~\cite{Karatsuba}. 

We also conjecture that
 the correction term $\theta_n$ is
\emph{completely} monotonic.

\subsection*{Acknowledgements}

We would like to thank Daniel Campos-Salas for some helpful
suggestions about the concavity. MBV acknowledges support for the
Vicerrector\'ia de Investigaci\'on of the University of Costa Rica.




\begin{thebibliography}{4}

\bibitem{HirschhornNewStir}
M. D. Hirschhorn,
``A new version of Stirling's formula'',
Mathl. Gazette \textbf{90} (2006), 286--291.

\bibitem{Karatsuba}
E. A. Karatsuba,
``On the asymptotic representation of the Euler gamma function by
Ramanujan''.
J. Comput. Appl. Math. \textbf{135} (2001), 225--240.

\bibitem{RamanLost}
S. Ramanujan,
\textit{The Lost Notebook and other Unpublished Papers},
S. Raghavan and S. S. Rangachari, eds.,
Narosa, New Delhi, 1987.

\bibitem{VillarinoCC}
M. B. Villarino, D. Campos-Salas and J. Carvajal-Rojas,
``On the monotonicity of the correction term in Ramanujan's factorial
approximation'',
Mathl. Gazette \textbf{97} (2013), to appear.

\end{thebibliography}
\end{document}